\documentclass[12pt,a4paper,draft]{article}
\usepackage[cp1251]{inputenc}
\usepackage[english,russian]{babel}
\usepackage{amssymb, amsmath, latexsym, amsfonts, amsthm}
\newtheorem*{lemma}{Lemma}
\newtheorem*{theorem}{Theorem}
\begin{document}
\renewcommand{\refname}{References}
\thispagestyle{empty}
\begin{center}
On weighted polynomial approximation 
\end{center}
\begin{center}
I.Kh. Musin
\end{center}
\begin{center}
{\it Institute of mathematics with Computer centre, Chernyshevskii str., 112, Ufa, 450077, Russia, musin\_ildar@mail.ru}
\end{center}

\renewcommand{\abstractname}{}
\begin{abstract}
{\sc Abstract}.  Let 
$\varPhi:{\mathbb R}^n \to [1, \infty)$  be a semi-continuous from below
function such that
$
\lim \limits_{x \to \infty}
\displaystyle \frac {\ln \varPhi(x)} {\Vert x \Vert} = +\infty. 
$
It is shown that polynomials are dense in $C_{\varPhi}({\mathbb R}^n)$.

\vspace {0.3cm}

\vspace {0.3cm}

\end{abstract}

\vspace{0,5cm}

{\bf 1. On a problem}. Let $\varPhi:{\mathbb R}^n \to [1, \infty)$ be such that   
$$
\lim \limits_{x \to \infty}\displaystyle \frac {\varPhi(x)} {\Vert x \Vert^m} = +\infty \ (m = 0, 1, 2, \ldots), 
$$
where $\Vert \cdot \Vert$ is the Euclidean norm in ${\mathbb R}^n$.
By $C_{\varPhi}({\mathbb R}^n)$ denote the normed  space of continuos functions in ${\mathbb R}^n$ such that $f(x) = o(\varPhi(x))$ as $x \to \infty$
with the norm
$$
p (f) = \sup_{x \in {\mathbb R}^n}
\displaystyle \frac
{\vert f(x) \vert}
{\varPhi(x)}.
$$

In 1924 S. Bernstein \cite {B} posed a problem to describe the weights $\varPhi$ such that polynomials are dense in $C_{\varPhi}({\mathbb R})$ and found a criterion of solvability of this problem in case when $\varPhi$ is the restriction to the real line of an even entire function with positive Maclaurin series coefficients. 

Bernstein’s problem was investigated by many prominent mathematicians. Independent solutions were obtained  by N.I. Ahiezer \cite {A}, S.N. Mergelyan \cite {M} and H. Pollard \cite {Pol} in 1953–1954 and by
L. de Branges \cite {Bra} in 1959. Information on the history of Bernstein’s problem can
be found in surveys by N.I. Ahiezer \cite {A}, S.N. Mergelyan \cite {M} (see also \cite {Lub}).

Bernstein's problem in several variables was considered in \cite {D-T}, \cite {S} for some particular cases.
It seems to be interesting to find easily verified conditions on $\varPhi$ for density of polynomials in $C_{\varPhi}({\mathbb R}^n)$.

{\bf 2. The main result}. Let ${\bf \Phi}$ be a family of semi-continuous from below
functions 
$\varPhi:{\mathbb R}^n \to [1, \infty)$  
such that
$
\lim \limits_{x \to \infty}
\displaystyle \frac {\ln \varPhi(x)} {\Vert x \Vert} = +\infty. 
$

\begin{theorem}
Let ${\varPhi \in \bf \Phi}$. Then polynomials are dense in $C_{\varPhi}({\mathbb R}^n)$.
\end{theorem}

{\bf 3. Some notations}. 
Denote by ${\cal B}[0, \infty)$ a family of semi-continuous from below functions $g: [0, \infty) \to {\mathbb R}$ satisfying the condition 
$\lim \limits_{x \to +\infty} \displaystyle \frac {g(x)} {x} = +\infty. 
$

For $g \in {\cal B}[0, \infty)$ let $g^*$ be its the Young conjugate:
$$
g^*(x) = \sup \limits_{y \ge 0}(x y - g(y)), \
x \ge 0.
$$

If $[0, \infty)^n \subseteq X \subset {\mathbb R}^n$ then for a function $g$ on $X$ denote by $u[e]$ the function on  $[0, \infty)$ defined by the rule: $g[e](x) = g(e^x), \ x \ge 0$.

For $r >0$ let
$\Pi_r=\{x \in {\mathbb R}^n: \vert x_j \vert < r, j =1, \ldots , n \}$.

{\bf 4. Auxiliary result}. 

\begin{lemma}
Let $g \in {\mathcal B}[0, \infty)$. 
Then
$$
(g[e])^*(x) + (g^*[e])^*(x) \le x \ln x - x, \
x > 0.
$$
$$
(g[e])^*(0) + (g^*[e])^*(0) \le 0.
$$
\end{lemma}

{\bf Proof}. Let $g \in {\mathcal B}[0, \infty)$. Then $g^* \in {\mathcal B}[0, \infty)$. Next, 
let $x \ge 0$. For each $\varepsilon > 0$ there are points $t \ge 0$ и $\xi \ge 0$ such that 
$$
(g[e])^*(x) < x t - g[e](t) + \varepsilon, 
$$
$$
(g^*[e])^*(x) < x  \xi  - g^*[e](\xi) + \varepsilon.
$$
Hence,
$$
(g[e])^*(x) + (g^*[e])^*(x) < x t - g[e](t) + x  \xi - \sup_{\eta \ge 0} (e^{\xi} \eta - g(\eta)) + 2  \varepsilon. 
$$
Thus, for any $\eta \ge 0$ 
$$
(g[e])^*(x) + (g^*[e])^*(x) < x t - g(e^t) + x  \xi - e^{\xi} \eta + g(\eta) + 
2 \varepsilon.
$$
Letting here $\eta = e^t$ we get
$$
(g[e])^*(x) + (g^*[e])^*(x) <  x t  + x  \xi - e^{\xi + t} + 2 \varepsilon. 
$$
Consequently,
$$
(g[e])^*(x) + (g^*[e])^*(x) <  \sup_{y \ge 0} (x y - e^y) + 2  \varepsilon \le \sup_{y \in {\mathbb R}} (x y - e^y) + 2  \varepsilon = x \ln x - x + 2  \varepsilon. 
$$
From this (since $\varepsilon > 0$ is arbitrary) we obtain that 
$$
(g[e])^*(x) + (g^*[e])^*(x) \le x \ln x - x, \ x > 0;
$$
$$
(g[e])^*(0) + (g^*[e])^*(0) \le 0. 
$$

{\bf 5. Proof of the Theorem}. Let $f \in C_{\varPhi}({\mathbb R}^n)$. 
Let us approximate $f$ by polynomials in $C_{\varPhi}({\mathbb R}^n)$.
There are three steps in the proof.

1. Choose a function $\chi \in C^{\infty}({\mathbb R})$  such that
supp  $ \chi \subseteq [-2,2],
\chi (x) = ~1 $ для $ x \in [-1,1]$, $ 0 \le  \chi (x) \le 1 $
$ \forall x \in {\mathbb R}$.
Put
$
\eta (x_1, x_2, \ldots , x_n) = \chi (x_1) \chi (x_2) \cdots \chi (x_n).
$

Let
$ f_{\nu}(x) = f(x)\eta (\frac x {\nu})$, \
${\nu} \in {\mathbb N}$,
$x \in {\mathbb R}^n $.
It is clear that $f_{\nu} \in C_{\varPhi}({\mathbb R}^n)$.
Since
$$
\sup \limits_{x \in {\mathbb R}^n} \displaystyle \frac
{\vert f_{\nu}(x) - f(x) \vert} {\varPhi(x)} \le
\sup \limits_{x \notin \Pi_{\nu}}
\displaystyle \frac {\vert f(x)
\vert}{\varPhi(x)} ,
$$
then $p(f_{\nu} - f) \to 0 $
as $ \nu \to \infty $.
Thus, the sequence $(f_{\nu})_{\nu=1}^{\infty}$ converges to $f$ in 
$C_{\varPhi}({\mathbb R}^n)$ as $ \nu \to \infty $.

2. Fix  $\nu \in {\mathbb N}$. 
Let $h$ be an entire function of exponential type at most 1 such that 
$h \in L_1({\mathbb R})$ and $h(x) \ge 0$ for $x \in  {\mathbb R}$. 
For example, we can take $h(z) = \displaystyle \frac  {\sin^2{ \frac z 2}} {z ^2}, \
z \in~{\mathbb C}$. 
Put
$H(z_1, z_2, \ldots , z_n) = h(z_1) h(z_2) \cdots h(z_n)$.
Using the Paley-Wiener theorem we can find a constant $K_H >0$ such that for each 
$\alpha \in {\mathbb Z_+^n}$
\begin{equation}
\vert
(D^{\alpha} H)(x) \vert \le K_H \ , \ x \in {\mathbb R}^n.
\end{equation}
Let
$ \int_{{\mathbb R}^n} H(x) \ d x = A $.
For $\lambda> 1$ let
$$
f_{\nu,\lambda}(x) =  \displaystyle \frac  {{\lambda}^n} A
\int_{{\mathbb R}^n} f_{\nu}(y)
H(\lambda (x-y)) \ d y, \ x \in {\mathbb R}^n.
$$
Obviously, $f_{\nu,\lambda} \in C({\mathbb R}^n)$ and $f_{\nu,\lambda}$ is bounded. 
Hence, $f_{\nu,\lambda} \in C_{\varPhi}({\mathbb R}^n)$.
Let us show that $f_{\nu,\lambda} \to f_{\nu} $ in $C_{\varPhi}({\mathbb R}^n)$ as
$\lambda \to +\infty $. 
Let $\varepsilon > 0$ be arbitrary. 
For $\lambda> 1$ let
$r(\lambda) = {\lambda}^{-\frac {2n}{2n+1}}$.
For any $ x \in {\mathbb R}^n$ we have that
$$
f_{\nu, \lambda}(x) - f_{\nu}(x)
=
\displaystyle \frac  {{\lambda}^n} A
\int_{{\mathbb R}^n}
(f_{\nu}(y) - f_{\nu}(x))
H(\lambda (x-y)) \ dy =
$$
$$
= \displaystyle \frac  {{\lambda}^n} A
\int_{\Vert y - x \Vert \le r(\lambda)}
(f_{\nu}(y) - f_{\nu}(x))
H(\lambda (x-y)) \ dy + 
$$
$$
+ \displaystyle \frac  {{\lambda}^n} A
\int_{\Vert y - x \Vert > r(\lambda)}
(f_{\nu}(y) - f_{\nu}(x))
H(\lambda (x-y)) \ dy .
$$
Denote the terms on the right-hand side of this equality by
$I_{1, \alpha}(x)$ and $I_{2, \alpha}(x)$, respectively.

Let us estimate $I_{1, \alpha}(x)$. In view of uniform continuity of function $f_{\nu}$ in ${\mathbb R}^n$ there is  a number $\delta = \delta(\varepsilon)> 0$ such that for any  
$t, u \in {\mathbb R}^n$ such that $\Vert t - u \Vert < \delta $ we have that 
$\vert f_{\nu}(t) - f_{\nu}(u) \vert < \varepsilon $. Then for all  
$\lambda > \lambda_1: =\delta^{-\frac {2n+1}{2n}}$ and all $x \in {\mathbb R}^n$ we get that
$$
\vert I_{1, \alpha}(x) \vert 
< \displaystyle \frac  {\varepsilon \lambda^n} {A} \int_{\Vert y - x \Vert \le r(\lambda)}
H(\lambda (x-y)) \ dy \le \displaystyle \frac  {\varepsilon \lambda^n} {A}
\int_{{\mathbb R}^n} 
H(\lambda (x-y)) \ d y = \varepsilon
$$
Let us estimate $I_{2, \alpha}(x)$. Let $M_{\nu} = \max \limits_{\Vert x \Vert \le 2 \nu} \vert f(x) \vert$.  
Then for any $x \in {\mathbb R}^n$ 
$$
\vert I_{2, \alpha}(x) \le 
\displaystyle \frac  {2  M_{\nu} {\lambda}^n} {A} 
\int_{\Vert y - x \Vert > r(\lambda)}
H(\lambda (x-y)) \ dy = 
\displaystyle \frac  {2  M_{\nu}} {A} \int_{\Vert t \Vert > {\lambda}^{\frac 13}}
H(t) \ dt .
$$
Since $H \in L_1({\mathbb R}^n)$ then there is a number $\lambda_2 > 0$ such that for all 
$\lambda > \lambda_2$ we have that $\vert I_{2, \alpha}(x) \vert < \varepsilon$ for any $x \in {\mathbb R}^n$. 
Thus, if $\lambda > \max (\lambda_1, \lambda_2)$ then 
$\vert f_{\nu, \lambda}(x) - f_{\nu}(x) \vert < 2 \varepsilon$ for any 
$x \in {\mathbb R}^n$. This means that
$p(f_{\nu,\lambda} - f_{\nu}) \le 2 \varepsilon $ for $\lambda > \max (\lambda_1, \lambda_2)$.
Thus, $f_{\nu,\lambda} \to f_{\nu}$ in $C_{\varPhi}({\mathbb R}^n)$ as
$\lambda \to +\infty $. 

3. Fix $\lambda>0, {\nu} \in {\mathbb N} $.  let us approximate
$f_{\nu,\lambda}$ by polynomials in $C_{\varPhi}({\mathbb R}^n)$.
For $N \in {\mathbb N}, x = (x_1, \ldots , x_n) \in {\mathbb R}^n$ let
$$
U_N(x) = H(0) + \displaystyle \sum \limits_{k=1}^{N}
\frac {\displaystyle \sum_{1\le i_1 \le n} \cdots
\sum_{1\le i_k \le n} \displaystyle \frac {{\partial}^k H}
{\partial x_{i_1} \cdots \partial x_{i_k}}(0)
x_{i_1} \cdots x_{i_k}} {k!} \ .
$$Since for $x \in {\mathbb R}^n$
$$
\vert H(x) - U_N(x) \vert  \le
\frac {\displaystyle \sum_{1\le i_1 \le n} \cdots
\sum_{1\le i_{N+1} \le n} 2 \sup_{\xi \in (0, x)}
\left \vert \displaystyle \frac {{\partial}^{N+1} H}
{\partial x_{i_1} \cdots \partial x_{i_{N+1}}}(\xi)
x_{i_1} \cdots x_{i_{N+1}} \right \vert} {(N+1)!} \ ,
$$
then using the inequality (1) we obtain that
\begin{equation}
\vert H(x) - U_N(x) \vert \le
\displaystyle \frac {2K_H n^{N+1} {\Vert x \Vert}^{N+1}}{(N+1)!} \ .
\end{equation}
Put
$$
V_N(x) =
\displaystyle \frac  {{\lambda}^n} {A}
\int_{\Pi_{2 \nu}}
f_{\nu}(y) U_{N}(\lambda (x-y)) \ d y , \ x \in {\mathbb R}^n.
$$
It is clear that $V_N$ is a polynomial of degree at most $N$. 
We claim that the sequence $(V_N)_{N=1}^{\infty}$ converges to
$ f_{\nu,\lambda} $ in $C_{\varPhi}({\mathbb R}^n)$ as $N \to \infty $.
For $x \in {\mathbb R}^n$ we have that
$$
f_{\nu,\lambda}(x) - V_{N}(x))
=
\displaystyle \frac {{\lambda}^n} {A}
\int_{\Pi_{2 \nu}}
(f_{\nu}(y)
(H({\lambda}(x-y)) - U_N(\lambda (x-y))) \ d y .
$$
From this using (2) we can find positive constants $C_1$ and $C_2$ (not depending on $N$) 
such that for all $N \in {\mathbb  N}$, 
$x \in {\mathbb R}^n $
$$
\vert
f_{\nu,\lambda}(x) -  V_{N}(x)
\vert
\le \displaystyle \frac {C_1 C_2^N (1+\Vert x \Vert)^{N+1}}
{(N+1)!} \ .
$$
Thus, for each
$N {\in {\mathbb  N}}$
\begin{equation}
p(f_{\nu,\lambda} - V_N) \le
\displaystyle \frac
{C_1 C_2^N}{(N+1)!}
\sup \limits_{x \in {\mathbb R}^n}
\displaystyle \frac {(1+\Vert x \Vert)^{N+1}} {\varPhi(x)} \ .
\end{equation}
Put $\varphi(x) = \ln \varPhi (x), \ x \in {\mathbb  R}^n$. 
For
$\sigma \in S^{n-1}: = \{x \in {\mathbb R}^n: \Vert x \Vert = 1 \}$ let 
$\varphi_{\sigma}(t) =
\varphi(\sigma t), \ t \ge 0.
$
Note that
$$
\sup \limits_{x \in {\mathbb R}^n}
\displaystyle \frac {(1+\Vert x \Vert)^{N+1}} {\varPhi(x)} = 
\exp (\sup \limits_{r \ge 0, \sigma \in S^{n-1}}((N+1)\ln (r + 1) - \varphi (r \sigma))) \le 
$$
$$
\le 
2^{N+1} \exp (\max (0, \sup \limits_{r \ge 1, \sigma \in S^{n-1}}((N+1) \ln r - 
\varphi_{\sigma} (r)))) = 
$$
$$
=
2^{N+1} \exp (\max (0, \sup \limits_{\sigma \in S^{n-1}}
(\sup \limits_{r \ge 1}((N+1) \ln r - 
\varphi_{\sigma}(r))))) = 
$$
$$
=
2^{N+1} \exp (\max (0, \sup \limits_{\sigma \in S^{n-1}}
(\varphi_{\sigma}[e])^*(N+1))).
$$
Next, for each $\sigma \in S^{n-1}$ \ $\varphi_{\sigma} \in {\cal B}[0, \infty)$. Hence, by Lemma
$$
(\varphi_{\sigma}[e])^*(N+1)) \le (N+1) \ln (N+1) - (N+1) - ((\varphi_{\sigma})^*[e])^*(N+1).
$$
Therefore,
$$
\sup \limits_{\sigma \in S^{n-1}} 
(\varphi_{\sigma}[e])^*(N+1)) \le (N+1) \ln (N+1) - (N+1) - 
\inf \limits_{\sigma \in S^{n-1}}((\varphi_{\sigma})^*[e])^*(N+1).
$$
Thus, 
$$
\sup \limits_{x \in {\mathbb R}^n}
\displaystyle \frac {(1+\Vert x \Vert)^{N+1}} {\varPhi(x)}
\le 
2^{N+1} 
\max \left(1, \frac {(N+1)^{N+1}} 
{e^{N+1} e^{ \inf \limits_{\sigma \in S^{n-1}}
((\varphi_{\sigma})^*[e])^*(N+1)}}\right).
$$
Returning to (3) for each $N {\in {\mathbb  N}}$ we get that
$$
p(f_{\nu,\lambda} - V_N) \le 
\displaystyle \frac
{C_1 C_2^N 2^{N+1} }{(N+1)!}
\max \left(1, \frac {(N+1)^{N+1}} 
{e^{N+1} e^{\inf \limits_{\sigma \in S^{n-1}}((\varphi_{\sigma})^*[e])^*(N+1)}}\right).
$$
From this taking into account that $m! \ge \frac {m^m} {e^m}$ 
for $m \in {\mathbb N}$ we have  that
\begin{equation}
p(f_{\nu,\lambda} - V_N) \le 
2C_1   
\displaystyle \max \left(\frac {(2C_2)^N} {(N+1)!}, \frac {(2C_2)^N}
{e^{\inf \limits_{\sigma \in S^{n-1}}((\varphi_{\sigma})^*[e])^*(N+1)}}\right). 
\end{equation} 
Note that uniformly for $\sigma \in S^{n-1}$ 
\begin{equation}
\lim_{\xi \rightarrow+\infty} \frac {(\varphi^*_{\sigma}[e])^*(\xi)} {\xi} = +\infty 
\end{equation} 
since for each $\sigma \in S^{n-1}$ 
$$
(\varphi^*_{\sigma}[e])^*(\xi) \ge \xi t - \varphi^*_{\sigma}[e](t), \ \xi \ge 0, t \ge 0,
$$
and
$$
\varphi^*_{\sigma}[e](t) = \sup_{r \ge 0}
(e^t r - \varphi_{\sigma}(r)) \le \sup_{{r \ge 0}, \atop {\sigma \in  S^{n-1}}}
(e^t r - \varphi(r \sigma)) = 
\sup \limits_{x \in {\mathbb R}^n} (e^t \Vert x \Vert - \varphi(x)).
$$
In view of (5) from (4) it follows that
$
p(f_{\nu,\lambda} - V_N) \to 0
$
as $N \to \infty $.
Thus, the function $f_{\nu,\lambda}$ is approximated by polynomials in
$C_{\varPhi}({\mathbb R}^n)$.

From the conclusions of all three above steps, one derives that each function $f \in C_{\varPhi}({\mathbb R}^n)$ can be approximated by polynomials in $C_{\varPhi}({\mathbb R}^n)$.
 
This work was supported by the grant from RFBR 15-01-01661.

\end{document}